\documentclass{gtart_a}
\pdfoutput=1
\usepackage{pinlabel}

\title{Dense embeddings of surface groups}

\author[Breuillard]{Emmanuel Breuillard}
\givenname{Emmanuel}
\surname{Breuillard}
\address{Universite de Lille\\UFR de Mathematiques\\
59655 Villeneuve d'Ascq\\FRANCE}
\email{emmanuel.breuillard@math.univ-lille1.fr}
\urladdr{}

\author[Gelander]{Tsachik Gelander}
\givenname{Tsachik}
\surname{Gelander}
\address{Mathematics Department\\ Yale University\\
10 Hillhouse ave\\New Haven CT 06511\\USA}
\email{tsachik.gelander@yale.edu}
\urladdr{}

\author[Souto]{Juan Souto}
\givenname{Juan}
\surname{Souto}
\address{Dept of Maths\\University of Chicago\\
5734 S. University Avenue\\Chicago, IL 60637\\USA}
\email{juan@math.uchicago.edu}
\urladdr{}

\author[Storm]{Peter Storm}
\givenname{Peter}
\surname{Storm}
\address{Stanford University\\Mathematics, Bldg. 380\\
450 Serra Mall\\Stanford, CA 94305\\USA}
\email{storm@math.stanford.edu}
\urladdr{}

\volumenumber{10}
\issuenumber{}
\publicationyear{2006}
\papernumber{34}
\lognumber{0713}
\startpage{1373}
\endpage{1389}

\doi{}
\MR{}
\Zbl{}


\received{10 February 2006}
\revised{3 August 2006}
\accepted{18 June 2006}
\proposed{Benson Farb}
\seconded{Jean-Pierre Otal, Walter Neumann}
\published{4 October 2006}
\publishedonline{4 October 2006}
\corresponding{}
\editor{}
\version{}

\arxivreference{math.GR/0602635}

\AtBeginDocument{\let\bar\wbar}
\makeop{Inn}
\makeop{Aut}

\makeatletter
\def\cnewtheorem#1[#2]#3{\newtheorem{#1}{#3}[section]
\expandafter\let\csname c@#1\endcsname\c@prop}

\newtheorem{prop}{Proposition}[section]
\cnewtheorem{thm}[prop]{Theorem}
\cnewtheorem{lem}[prop]{Lemma}
\cnewtheorem{cor}[prop]{Corollary}
\cnewtheorem{conj}[prop]{Conjecture}
\cnewtheorem{theorem}[prop]{Theorem}
\cnewtheorem{proposition}[prop]{Proposition}
\cnewtheorem{lemma}[prop]{Lemma}
\cnewtheorem{clm}[prop]{claim}

\theoremstyle{definition}
\newtheorem*{Ack}{Acknowledgments}
\cnewtheorem{nota}[prop]{Notations}
\cnewtheorem{defn}[prop]{Definition}
\newtheorem*{rem}{Remark}
\cnewtheorem{exam}[prop]{Example}
\cnewtheorem{prob}[prop]{Problem}
\cnewtheorem{ques}[prop]{Question}
\cnewtheorem{ob}[prop]{Observation}
\cnewtheorem{remark}[prop]{Remark}
\cnewtheorem{definition}[prop]{Definition}
\cnewtheorem{conjecture}[prop]{Conjecture}

\makeatother  

\DeclareMathOperator{\Hom}{Hom}

\newcommand{\gC}{\Gamma}
\newcommand{\gc}{\gamma}

\DeclareMathOperator{\Ker}{Ker}

\begin{document}

\begin{abstract} 
We discuss dense embeddings of surface groups and fully residually free 
groups in topological groups. We show that a compact topological group 
contains a nonabelian dense free group of finite rank 
if and only if it contains a dense surface group. 
Also, we obtain a characterization of those Lie 
groups which admit a dense faithfully embedded surface group. Similarly, 
we show that any connected semisimple Lie group contains a dense copy of 
any fully residually free group.
\end{abstract}

\maketitle

\section{Introduction}

Given a locally compact topological group $G$ and an abstract group
$\Gamma $, it is natural to ask whether $\Gamma $ can be embedded
densely in $G$.  More generally, for a given $G$ one would like to
understand its dense subgroups, and for a given $\Gamma $ one would
like to know its possible completions $G$ which are topological
groups. These questions are more accessible in the case where $G$ is a
finite dimensional analytic Lie group over a local field. While
discrete subgroups of Lie groups have been thoroughly studied for the
last fifty years, very little is known about nondiscrete, and in
particular, dense\footnote{Note that when $G$ is a connected simple
Lie group, a generic subgroup with sufficiently many generators is
either discrete or dense} subgroups of Lie groups. A dense embedding
of $\Gamma $ in $G$ may yield interesting data on $\Gamma $, $G$,
and the spaces on which they act (Margulis \cite{mar}, Sullivan \cite{sul},
Gelander and \.{Z}uk \cite{GZ}, Lubotzky and Weiss \cite{LW},
Breuillard and Gelander \cite{BG1,BG2}, Abert and Glasner \cite{AG}).

By a \textit{surface group}, we mean the fundamental group of a closed
oriented surface of genus at least $2$. By a \textit{free group}, we mean a
nonabelian free group on at least two generators. We obtain various results, all of
which are proved by continuously deforming a given representation to a
faithful one.

Our first result states that free groups and surface groups have the same
compactifications within the category of topological groups.

\begin{thm}
\label{comp}\label{compact}Let $G$ be a compact group. Then the following
two assertions are equivalent:

\begin{itemize}
\item  $G$ contains a dense free subgroup of finite rank.

\item  $G$ contains a dense surface group.
\end{itemize}
\end{thm}

As a corollary we obtain that a compact group contains a surface group if and only if it contains a free group. It is sometimes fairly 
easy to verify that a given compact group contains a free subgroup by means of probabilistic methods. However, we do not know a simple 
characterization of the compact groups containing a dense free subgroup of finite rank. It was shown in \cite{BG2} that the profinite 
completion $\widehat{\Gamma }$ of a finitely generated linear group $\Gamma$ contains a dense free subgroup of finite rank if and only 
if it is not virtually solvable (ie, contains no solvable subgroup of finite index). However, there are examples of topologically 
finitely generated profinite groups that satisfy a nontrivial group law (hence admit no free subgroups) although they are not virtually 
solvable (de Cornulier and Mann \cite{Mann}). On the other hand, it is possible to verify that any connected second countable 
nonabelian compact group contains a dense free subgroup of rank 2 (see \fullref{fscg}). Hence any such group contains also a 
dense surface group.

The method used to prove our \fullref{compact} above can be pushed a little further to get a result that holds for an arbitrary locally compact group: 

\begin{thm}
\label{main}Let $G$ be a locally compact group. Suppose that $G$ contains a
nondiscrete free subgroup $F$ of finite rank $r>1$. Then $G$ has a subgroup
$\Gamma $ containing $F$ such that $\Gamma $ is isomorphic to a surface
group (of genus $2r$). In particular, if $G$ has a dense free subgroup of
finite rank, then it has a dense surface group.
\end{thm}

\begin{rem}
As a corollary of \fullref{main} we obtain an elementary proof of a
result from Gelander and Glasner \cite{GG}, that surface groups are
primitive, ie, admit faithful primitive permutation
representations. Indeed, let $\Gamma $ be a surface group and embed
$\Gamma $ densely in $\text{PSL}_{2}({\mathbb{Q}}_{p})$. Then
$\Delta =\Gamma \cap \text{PSL}_{2}({\mathbb{Z}}_{p})$ is a maximal
subgroup of $\Gamma $ which contains no nontrivial normal subgroup of
$\Gamma $, and the action of $\Gamma $ on $\Gamma /\Delta $ is
primitive and faithful.
\end{rem}

When $G$ is a (nondiscrete) real Lie group with a countable number of
connected components, then $G$
contains a finitely generated dense free group if and only if the connected
component of the identity $G^\circ$ is not solvable and $G/G^\circ$ is
finitely generated \cite{BG1,BG2}.

We thus obtain:

\begin{cor}
\label{charac}Let $G$ be a nondiscrete real Lie group. Then the following
are equivalent:

\begin{itemize}
\item  $G$ contains a finitely generated dense free subgroup.

\item  $G$ contains a dense surface group.

\item  $G^{\circ }$ is not solvable and $G/G^{\circ }$ is finitely generated.
\end{itemize}
\end{cor}

One key property of surface groups which motivated this research is 
the fact that they are fully residually free. A finitely 
generated group $\Gamma $ is \emph{fully residually free} if for 
every finite set $K \subset \Gamma \setminus \{1\}$ there is a 
homomorphism $\phi \co \Gamma \to {F}$ onto a free group $F$ with $K\cap\Ker(\phi)=\emptyset$. In other words, $\Gamma$ is fully residually 
free if any finite set can be separated through a surjective map onto 
a free group. 

The class of fully residually free groups also appears in 
the work of Sela \cite{sela}, where he shows that it 
coincides with his notion of limit groups. The fact that surface 
groups are fully residually free is due to Baumslag \cite{Baum}. A 
group $\Gamma $ is called \textit{$d$-fully residually free} if any 
finite set can be separated trough a surjection on ${F}_{d}$. Note 
that if $\Gamma $ is $d$-fully residually free then it is also 
$k$-fully residually free for any $2\leq k < d$ (see \fullref{baum} 
below).

For general fully residually free groups we prove the following:

\begin{thm} \label{frf} 
Let $G$ be a connected nonsolvable Lie group. Then there is a number $d=d(G)<\dim (G)$ (defined in \fullref{labelled section}) such that: if $\Gamma $ is a finitely generated $d$-fully residually free group, then there is a dense embedding 
$\Gamma \hookrightarrow G$. 
\end{thm}

When $G$ is topologically perfect, ie, does not surject onto the circle,
then we can take $d$ to be the minimal number of generators for the Lie algebra of $%
G $. Since any semisimple Lie algebra is generated by $2$ elements, we
obtain:

\begin{thm}
\label{semisimple} Any connected semisimple Lie group contains a dense copy
of any finitely generated nonabelian fully residually free group.
\end{thm}

\begin{rem}
A group $\Gamma $ is \emph{residually free} if for every $\gamma \in \Gamma
\setminus \{1\}$ there is $\phi \co \Gamma \to {F}_{d}$ with $\gamma \notin %
\Ker(\phi)$. For example, if $\Gamma $ is a surface group then $\Gamma
\times \Gamma $ is residually free. Since $\text{PSL}_{2}({\mathbb{C}})$
does not have subgroups isomorphic to $\Gamma \times \Gamma $ we observe
that in \fullref{frf} the condition ``$\Gamma $ is fully residually
free'' cannot be weakened to ``$\Gamma $ is residually free''.
\end{rem}

Let us end this introduction by remarking that all the results obtained in
this paper are concerned with the existence of subgroups with certain
desired properties. However we do not obtain concrete examples. In general,
this problem seems much more difficult.

\begin{Ack}
E Breuillard acknowleges support from the Centre National de la
Recherche Scientifique and from the Institute for Advanced Study.  T
Gelander has received support from NSF grant DMS-0404557 and BSF grant
2004010.  P Storm has received support from an NSF Postdoctoral
Research Fellowship.
\end{Ack}




\section{Eventually faithful homomorphisms and a lemma of Baumslag}

Let $(\rho _{n})_{n\geq 0}$ be a sequence of homomorphisms from a group $H$
to a group $G$. We say that $(\rho _{n})_{n\geq 0}$ is \textit{eventually
faithful} if for every $h\in H\backslash \{1\}$ there exists an integer $%
n_{0}=n_{0}(h)$ such that $\rho _{n}(h)\neq 1$ for all $n\geq n_{0}.$

Since any finitely generated free group can be embedded into $F_{2}$, the
free group on two generators, it follows that a finitely generated group is
fully residually free if and only if it admits an eventually faithful
sequence of homomorphisms to $F_{2}.$

Let us recall the following lemma of Baumslag \cite{Baum}.

\begin{lemma}[Baumslag]\label{baum} Let $u,a_{1},...,a_{k}$ be elements of a free
group $F$. Assume that $u$ does not commute with any of the $a_{i}$'s. Then
there exists $n_{0}\geq 0$ such that for all integers $n_{1},...,n_{k}$ with
$|n_{i}|\geq n_{0}$ we have
\begin{equation*}
u^{n_{1}}a_{1}u^{n_{2}}a_{2}\cdot ...\cdot u^{n_{k}}a_{k}\neq 1
\end{equation*}
\end{lemma}

This lemma has a few corollaries. The first one proves that surface groups are fully residually free.

\begin{cor}
\label{Juan}Let $\Gamma =\Gamma _{2r}$ be the fundamental group of an
orientable surface of genus $2r$ ($r\geq 1$). Let us write a presentation of
$\Gamma $ as
\begin{equation}
\Gamma =\left\langle a_{i},a_{i}^{\prime },b_{i},b_{i}^{\prime },1\leq i\leq
r \;\; | \;\; [a_{1},a_{1}^{\prime }]\cdot ...\cdot [a_{r},a_{r}^{\prime
}]\cdot [b_{r}^{\prime },b_{r}]\cdot ...\cdot [b_{1}^{\prime
},b_{1}]=1\right\rangle   \label{rep}
\end{equation}
Now consider the automorphism $\sigma $ of $\Gamma $ that leaves the $a_{i}$%
's and $a_{i}^{\prime }\,$'s fixed while sending every $b_{i}$ to $\gamma
b_{i}\gamma ^{-1}$ and every $b_{i}^{\prime }$ to $\gamma b_{i}^{\prime
}\gamma ^{-1}$, where $\gamma =\ [a_{1},a_{1}^{\prime }]\cdot ...\cdot
[a_{r},a_{r}^{\prime }].$ Finally let $f$ be the surjective homomorphism
from $\Gamma $ to the free group $F_{2r}$ with free generators $%
x_{1},...,x_{r}$ and $x_{1}^{\prime },...,x_{r}^{\prime }$ defined by $%
f(a_{i})=f(b_{i})=x_{i},$ $f(a_{i}^{\prime })=f(b_{i}^{\prime
})=x_{i}^{\prime }$.

Then the sequence of maps $(f\circ \sigma ^{n})_{n\geq 0}$ is eventually
faithful.
\end{cor}

\begin{figure}[ht!]
\begin{center}
\labellist\small\hair3pt
\pinlabel {$a_1$} [l] at 152 213
\pinlabel {$a'_1$} [bl] <0pt,-2pt> at 139 238
\pinlabel* {$a_1^{-1}$} [bl] at 120 255
\pinlabel {${a'}_1^{-1}$} [b] <2pt,0pt> at 97 263
\pinlabel {$a_2$} [b] <2pt,0pt> at 72 265
\pinlabel {$a'_2$} [br] at 49 255
\pinlabel* {$a_2^{-1}$} [br] at 31 238
\pinlabel* {${a'}_2^{-1}$} [r] at 20 213
\pinlabel {$b'_2$} [r] at 21 185
\pinlabel* {$b_2$} [tr] at 32 161
\pinlabel* {${b'}_2^{-1}$} [tr] at 52 141
\pinlabel {$b_2^{-1}$} [t] at 75 132
\pinlabel {$b'_1$} [t] at 101 131
\pinlabel {$b_1$} [tl] at 126 140
\pinlabel* {${b'}_1^{-1}$} [tl] at 144 159
\pinlabel {$b_1^{-1}$} [l] at 154 184
\pinlabel {$\gamma$} [b] at 87 198
\pinlabel {$f$} [l] at 88 100
\pinlabel {$x_1$} [l] at 154 16
\pinlabel {$x'_1$} [bl] <-3pt,0pt> at 141 42
\pinlabel {$x_1^{-1}$} [bl] at 121 58
\pinlabel {${x'}_1^{-1}$} [b] <4pt,0pt> at 98 68
\pinlabel {$x_2$} [b] at 75 68
\pinlabel {$x'_2$} [br] <2pt,0pt> at 50 58
\pinlabel* {$x_2^{-1}$} [br] at 32 41
\pinlabel* {${x'}_2^{-1}$} [r] at 21 14
\endlabellist
\includegraphics{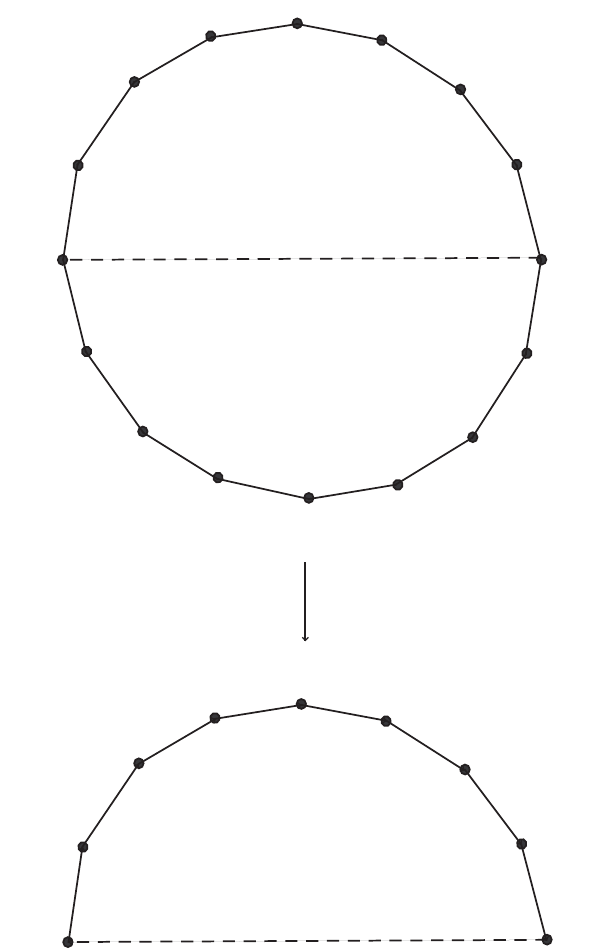}
\end{center}
\caption{All curves are oriented counterclockwise. Fold the genus $4$
surface across $\gamma$ to obtain $f$.}
\label{cool picture by Pete}
\end{figure}

The maps $\sigma $ and $f$ have the following simple topological
interpretation. In the above classical representation of $\Gamma$ as the
fundamental group of a surface, the relation gives the gluing instructions
for forming a genus $2r$ surface from a $4r$-gon in the plane. This is
demonstrated in \fullref{cool picture by Pete} for the case $r=2$. The
element $\gamma$ correponds to the closed curve separating the surface into
two equal parts. The map $\sigma $ corresponds to a Dehn twist around $%
\gamma $. The map $f$ is obtained by reflecting the surface across the
separating curve $\gamma$. The image of this reflection is a surface of
genus $r$ with one boundary component, whose fundamental group is freely
generated by $x_1, ..., x_r, x_1^{\prime}, ..., x_r^{\prime}$.

\proof
Let $g\in \Gamma \backslash \{1\}.$ The element $g$ can be written in the
form $g=w_{1}(a_{i},a_{i}^{\prime })\cdot w_{2}(b_{i},b_{i}^{\prime })\cdot
..\cdot w_{2p-1}(a_{i},a_{i}^{\prime })\cdot w_{2p}(b_{i},b_{i}^{\prime })$
where each $w_{i}$ is a reduced word in $2r$ letters and the first and the
last $w_{i}$ may be trivial. Up to modifying the odd $w_{2j-1}$, we may
assume that each even $w_{2j}$ ($1\leq j\leq p$) is such that $%
w_{2j}(b_{i},b_{i}^{\prime })$ is not a power of $\gamma $. Note that the
centralizer of $\gamma $ in $\Gamma $ is the cyclic group generated by $%
\gamma $. By regrouping several $w_{j}$'s into a longer word if necessary,
unless $g$ itself is a power of $\gamma $, we may also assume that 
$w_{2j-1}(a_{i},a_{i}^{\prime })$ is not a power of $ \gamma $. 
Let $\overline{\gamma }$ be the image of $\gamma $ under $f$. We
have $f\circ \sigma ^{n}(g)=w_{1}\cdot \overline{\gamma }^{n}w_{2}\overline{%
\gamma }^{-n}\cdot ..\cdot w_{2p-1}\cdot \overline{\gamma }^{n}w_{2p}%
\overline{\gamma }^{-n}$ where each $w_{j}=w_{j}(x_{i},x_{i}^{\prime }).$
Since $\gamma $ does not commute with any of the $w_{j}$, \fullref{baum}
implies that $f\circ \gamma ^{n}$ is eventually faithful. \endproof

The next two corollaries are very simple applications of \fullref{baum},
and are only recorded here for further use.

\begin{cor}
\label{enlargeF}Let $F$ be a free group of rank $n+1$ with free generators $%
x_{1},...,x_{n+1}.$ Let $F^{-}$ be the subgroup generated by $%
x_{1},...,x_{n}.$ Suppose $a$ and $b$ are noncommuting elements in $F^{-}.$
Consider the automorphism $\sigma $ of $F$ defined by $\sigma (x_{i})=x_{i}$
if $i\leq n$ and $\sigma (x_{n+1})=bx_{n+1}b^{-1}.$ Let $f$ be the
homomorphism of $F$ into $F^{-}$ that sends each $x_{i}$ to itself for
 $1\leq i\leq n,$
and $x_{n+1}$ to $a$. Then the sequence of homomorphisms $(f\circ
\sigma ^{n})_{n\geq 0}$ is eventually faithful.
\end{cor}

\begin{rem}
It follows from the definition that \fullref{baum} remains true when the
free group $F$ is replaced by any nonabelian fully residually free group,
and in particular by a surface group.
\end{rem}

\begin{cor}
\label{conv}Let $\Gamma _{r}=\left\langle a_{i},a_{i}^{\prime },1\leq i\leq r%
\;\; | \;\; [a_{1},a_{1}^{\prime }]\cdot ...\cdot [a_{r},a_{r}^{\prime
}]=1\right\rangle $ be a presentation of a surface group of genus $r$. Let $F
$ be a free group of rank $2r$ generated by $x_{i},x_{i}^{\prime }$ for $%
1\leq i\leq r$. For each integer $n\geq 0$ consider the homomorphism $\rho
_{n}\co F\rightarrow \Gamma _{r}$ given by $\rho _{n}(x_{i})=a_{i}$ for $i\geq 2
$, $\rho _{n}(x_{i}^{\prime })=a_{i}^{\prime }$ for $i\geq 1$ and $\rho
_{n}(x_{1})=a_{1}\cdot (a_{1}^{\prime })^{n}.$ Then the sequence $(\rho
_{n})_{n\geq 1}$ is eventually faithful.
\end{cor}




\section[Proof of \ref{comp}]{Proof of \fullref{comp}}

Here we give a proof of \fullref{compact}. Let $G$ be a compact group
containing a dense free group $F$ on $r$ generators. We are going to show
that $G$ contains a surface group containing $F$. In order to do so, we
first make sure that $G$ contains a dense free group on an even number of
generators. This is done, if $r$ is odd, by enlarging $F$ in the following
way. Let $a_{1},...,a_{r}$ be generators of $F$ and fix $a$ and $b$ in $F$
two noncommuting elements. Then let $B$ be the closure in $G$ of the cyclic
group generated by $b$. Let $F^{+}$ be the abstract free group on $r+1$
generators $y_{1},...,y_{r+1}.$ To every $\beta \in B,$ we associate the
homomorphism $\rho _{\beta }\co F^{+}\rightarrow G$ which sends each $%
y_{i}$ to $a_{i}$ when $1\leq i\leq r$, and sends $y_{r+1}$ to $\beta a\beta
^{-1}.$ From \fullref{enlargeF}, we know that the sequence of
homomorphisms $(\rho _{b^{n}})_{n\geq 1}$ is eventually faithful. Let $w\in
F^{+}\backslash \{1\}$ and consider the set $O_{w}=\{\beta \in B \; : \; \rho
_{\beta }(w)\neq 1\}.$ Clearly $O_{w}$ is open in $B$. It is also dense
because the set $\{b^{n} \; : \; n\geq n_{0}\}$
 is dense in $B$ for any $n_{0}\geq 0$. Baire's theorem implies that 
$O=\cap _{w\in F^{+}\backslash
\{1\}}O_{w}$ is dense in $B$, and is in particular nonempty. Let $\beta _{0}\in
O.$ The homomorphism $\rho _{\beta _{0}}$ is faithful and $\rho _{\beta
_{0}}(F^{+})$ is a dense free subgroup of $G$ of rank $r+1$.

We may therefore assume that $r=2k.$ Let $x_{1},x_{1}^{\prime
},...,x_{k},x_{k}^{\prime }$ be the $2k$ free generators of $F.$ Set $\gamma
=[x_{1},x_{1}^{\prime }]\cdot ...\cdot [x_{k},x_{k}^{\prime }].$ Let $K$ be
the closure in $G$ of the cyclic group generated by $\gamma $. Keeping the
same notation as in \fullref{Juan} for the presentation of the
surface group $\Gamma _{2r}$, we define for every $\alpha \in K$ a
homomorphism $\sigma _{\alpha }\co \Gamma _{2r}\rightarrow G$ by sending $a_{i}$
to $x_{i}$, $a_{i}^{\prime }$ to $x_{i}^{\prime }$, $b_{i}$ to $\alpha
x_{i}\alpha ^{-1}$ and $b_{i}^{\prime }$ to $\alpha x_{i}^{\prime }\alpha
^{-1}.$ From \fullref{Juan}, we know that the sequence of
homomorphisms $(\sigma _{\gamma ^{n}})_{n\geq 1}$ is eventually faithful. As
above, let $w\in \Gamma _{2r}\backslash \{1\}$ and consider the set $%
U_{w}=\{\alpha \in K \; | \; \sigma _{\alpha }(w)\neq 1\}.$ Clearly $U_{w}$ is open
in $K$. It is also dense, because for any $n_{0}\geq 0$ the set $\{\gamma
^{n} \; : \; n\geq n_{0}\}$ is dense in $K$. Applying Baire's theorem, we obtain
that $U:=\cap _{w\in F^{+}\backslash \{1\}}U_{w}$ is dense in $K$ and in
particular nonempty. Let $\alpha _{0}\in U.$ Note that for every $\alpha
\in K$ the image $\sigma _{\alpha }(\Gamma _{2r})$ is dense in $G$ because
it contains $F$ as a subgroup. The homomorphism $\sigma _{\alpha _{0}}$ is
faithful and $\sigma _{\alpha _{0}}(\Gamma _{2r})$ is a surface group
densely embedded in $G.$

We now pass to the converse statement. Let $\Gamma _{r}$ be a dense surface
group of genus $r$ in $G.$ Keep the notation of \fullref{conv} and
let $A$ be the closure in $G$ of the cyclic group generated by $%
a_{1}^{\prime }.$ For every $\omega \in A$ let $\pi _{\omega
}\co F_{2r}\rightarrow G$ be the homomorphism that sends $x_{1}$ to $%
a_{1}\omega $, while for $i\geq 2,$ $x_{i}$ is sent to $a_{i}$ and for $%
i\geq 1,$ $x_{i}^{\prime }$ is sent to $a_{i}^{\prime }.$ According to
\fullref{conv}, the sequence $(\pi _{(a_{1}^{\prime })^{n}})_{n\geq 1}$
is eventually faithful. A Baire argument similar to the one above show that $%
\pi _{\omega _{0}}$ is faithful for some $\omega _{0}\in A.$ It remains to
check that $\pi _{\omega _{0}}(F_{2r})$ is dense in $G.$ This is clear
because it contains all $a_{i}$'s for $i\geq 2$ and $a_{i}^{\prime }$ 's for
$i\geq 1$. In particular $\pi _{\omega _{0}} (F_{2r})$ contains $a_{1}^{\prime }$,
hence the closure of $\pi _{\omega _{0}}(F_{2r})$ must contain $\omega _{0}.$
Therefore the closure of $\pi _{\omega _{0}}(F_{2r})$ must contain $a_{1}$,
implying it is all of $G$. This completes the proof of \fullref{comp}.




\section[Proof of \ref{main}]{Proof of \fullref{main}}
Let us make the obvious remark that there are (non--locally compact) Hausdorff nondiscrete topological groups where the property of Theorem \ref {main} does not hold. For instance consider the free group with the induced topology coming from a dense embedding inside a compact Lie group.

Let $x_{1},...,x_{r}$ be the $r$ free generators of the nondiscrete free
subgroup $F$.  As in the compact case, we are first going to enlarge the free subgroup $F$
(of rank $r$) to a bigger free subgroup of rank $2r$ by adding $r$ free
generators, then deform that free subgroup into a surface group. 

By the
structure theory of locally compact groups (Van Dantzig's theorem, see Montgomery
and Zippin \cite{MZ}) $G$ has an open subgroup $H'$ containing the connected component of the identity $G^{\circ }$ in such a way that $H'/G^{\circ }$ is compact. Moreover there is a normal compact subgroup $K$ of $H'$ such that $H'/K$ is a Lie group,
 and $K$ can be chosen sufficiently small for the finite set of conjugates $\{K_i = x_i K x_i^{-1} \}_{i=1}^r$ to be contained in the open finite intersection $\cap_{i=1}^r x_i H' x_i^{-1}$ \cite[Theorem 4.6]{MZ}.  
This second assumption is used only to know that the set $K K_1 \ldots K_r$ is a subgroup of $G$.  Up to replacing $H'$ by a smaller open subgroup $H \le G$ we can assume that $H/K$ is connected.  

Let $U$ and $V$ be sufficiently small neighborhoods of the identity in $H$
so that $x_{i} u x_{i}^{-1} u^{-1} \in V$ for any $u\in U$
and $i=1,...,r$, and so that the projection of any element in $V^{r}$ lies in
 a $1$-parameter
subgroup of $H/K$. We are going to find elements $x_{1}^{\prime
},...,x_{r}^{\prime }$ in $U$ which,
together with the $x_{i}$'s, form $2r$ free generators of a free
subgroup of $G$.

For this purpose, pick two noncommuting elements $a$ and $b$ in $F$ that are in $U$ modulo $K$.
This is always possible because $F$ is not discrete.
The proof will have two cases. Case (I) is when $F\cap K\neq \{1\}$. Case 
(II) is when $F\cap K=\{1\}$. In case (I) we can clearly assume that $a$
and $b$ belong to $K$. (Pick an element in $F\cap K$ and some suitable
conjugate of it.) Suppose that $x_{1}^{\prime },...,x_{j}^{\prime }$ have
been constructed. Define $F_{r+j+1}$ to be an abstract free group on $r+j+1$
generators $y_{1},...,y_{r},y_{1}^{\prime },...,y_{j+1}^{\prime }$, and let
us find $x_{j+1}^{\prime }.$  We will handle the two cases separately.

Assume first we are in case (I).  Let $B \le K$ be the closure of the cyclic group generated by $b$.  For $\beta \in B$ let $\rho_\beta $ be the homomorphism sending each $y_i $ to $x_i $, each $y'_i$ to $x'_i $ for $i \le j$, and $y'_{j+1}$ to $\beta a \beta^{-1} $.  \fullref{enlargeF} and Baire's theorem ensure that the subset of those $\beta $ for which $ \rho_\beta $ is faithful is Baire dense in $B$, and hence nonempty.  Fix such an element $\beta_0 \in B$.  The desired new generator is $x'_{j+1} = \rho_{\beta_0} ( y'_{j+1})$.

Now assume we are in case (II).  Morally, we repeat the argument of case (I), but the details differ.  By induction we may assume the free group $ \langle x_1$, $\ldots$, $x_r$, $ x'_1$, $ \ldots$ , $x'_j \rangle$ intersects $K$ trivially.  Let $B$ be the $1$-parameter subgroup of $H/K$ containing the coset $bK$.  For $\beta \in B$ let $\rho_\beta $ be the map to $H/K$ sending each $y_i$ to $x_i K$, each $y'_i $ to $x'_i K $ for $i \le j$, and $y'_{j+1}$ to $\beta \overline{a} \beta^{-1} $, where $ \overline{a} = aK $.  This yields a one parameter family of representations from $\langle y_1, \ldots, y_r, y'_1, \ldots, y'_{j+1} \rangle $ to $H/K$.  For any word $w$ in $\langle y_1, \ldots, y_r, y'_1, \ldots, y'_{j+1} \rangle $, the set $\{ \beta \in B \; : \; \rho_\beta (w) \not= 1 \}$ is open because the map $\text{ev}_w \co  \beta \mapsto \rho_\beta (w)$ is continuous.  By \fullref{enlargeF} this continuous map $\text{ev}_w$ is not constant.  Now we use the fact that $B$ and $H/K$ are Lie groups, which implies they are real analytic manifolds and $\text{ev}_w$ is real analytic.  Therefore the closed set $\text{ev}_w^{-1} (1) $ is nowhere dense.  By Baire's theorem the subset of $\beta \in B$ for which $\rho_\beta $ is faithful is nowhere dense in $B$.  Fix such a $\beta_0 \in B$ sufficiently near the identity for there to be an element $x'_{j+1}$ contained in the intersection $\rho_{\beta_0}^{-1} (y'_{j+1}) \cap U  \subset G$.  This choice of $x'_{j+1}$ completes the induction in case (II).


Continuing the argument in both cases, call $F^{\prime }$ the new free subgroup on $2r$ generators. Note that $F\leq F^{\prime }$,
 and in case (II) we have $F^{\prime }\cap K=\{1\}.$
Consider the product of commutators $\gamma =[x_{1},x_{1}^{\prime
}]\cdot ...\cdot [x_{r},x_{r}^{\prime }]$. Let $\Gamma_{2r}$ be a surface
group given with the presentation written above in $\eqref{rep}.$ Consider
the centralizer $Z_{G}(\gamma )$ of $\gamma $ in $G$. Given an element $%
\alpha $ in $Z_{G}(\gamma )$ we can define a representation $\rho _{\alpha
}\co \Gamma _{2r}\rightarrow G$ by setting $\rho _{\alpha }(a_{i})=x_{i},$ $%
\rho _{\alpha }(a_{i}^{\prime })=x_{i}^{\prime }$, $\rho _{\alpha
}(b_{i})=\alpha x_{i}\alpha ^{-1}$, and $\rho _{\alpha }(b_{i}^{\prime
})=\alpha x_{i}^{\prime }\alpha ^{-1}.$ \fullref{Juan} shows that the
sequence $(\rho _{\gamma ^{n}})_{n\geq 1}$ is eventually faithful. We will
make use of the following lemmas:

\begin{lemma}
\label{lift}Let $H$ be a locally compact group and $K$ a compact normal
subgroup such that $H/K$ is a Lie group. Let $\{x(t)\}_{t}$ be a $1$%
-parameter subgroup in $H/K$. Then $\{x(t)\}_{t}$ can be lifted to a $1$%
-parameter subgroup $\{\widetilde{x}(t)\}_{t}$ in $H$ such that $\pi (%
\widetilde{x}(t))=x(t)$ where $\pi \co H\rightarrow H/K$ is the quotient map.
\end{lemma}

\proof
See the end of Section 4.7 of \cite{MZ}. \endproof

\begin{lemma}
\label{centralizer} Let $H$ be a locally compact group and $K$ a
compact normal subgroup such that $H/K$ is connected. Then $H=Z_{H}(K)K$
where $Z_{H}(K)$ is the centralizer of $K$ in $H$.
\end{lemma}

\proof
It follows from \fullref{lift} that $H=H^\circ K$.

Let $\rho \co H\rightarrow \Aut(K)$ be a map that sends $h\in H$ to the automorphism 
$\text{i}(h) $ of $K$ given by the conjugation by $h$. Then $\ker \rho =Z_{H}(K).$
We need to show that $\rho (H)=\rho (K)$ and for this it is clearly enough
to prove that $\rho (H^{\circ })$ is contained in $\Inn(K)$, the group of
inner automorphisms of $K$. This is a consequence of the following lemma:

\begin{lemma}
\label{centralizer2} Let $K$ be a compact group. Then the connected
component of the identity of $\Aut(K)$ is contained in $\Inn(K).$
\end{lemma}

\proof First assume that $K$ is a Lie group. Then $K\simeq D\times T\times S$
where $S$ is semisimple, $T$ a torus and $D$ is a finite group. As is
well-known, $\Inn(S)$ has finite index in $\Aut(S)$, and $\Aut(T)$ is discrete.
It follows easily that $\Aut(K)^\circ\leq \Inn(S)$. Now we pass to the general
case.

According the Peter-Weyl theorem, $K$ has a descending chain of compact
normal subgroups $C_{1}\supset C_{2}\supset C_{3}\supset \ldots $ such that
any open neighborhood of the identity contains all but a finite number of
the subgroups $\{C_{i}\}$, and the quotient $K/C_{i}$ is always a Lie group.
By pulling back a small identity neighborhood from $K/C_{i}$ to $K$ we
obtain an open set $U_{i}\subset K$ containing $C_{i}$ such that any
subgroup of $K$ inside $U_{i}$ is in fact contained in $C_{i}$. Therefore by
connectivity every automorphism in $\Aut(K)^{\circ }$ preserves $C_{i}$. This
yields a map
\begin{equation*}
\Aut(K)^{\circ }\rightarrow \Aut(K/C_{i})^{\circ }\le \Inn(K/C_{i}).
\end{equation*}

Let $\phi \in \Aut(K)^\circ$. For each $i$ pick an element $h_i \in K$ such
that $\phi (g) C_i = h_i g h_i^{-1} C_i$ for all $g \in K$.
Since the subgroups $\{ C_i \}$ become arbitrarily small it follows that $%
\text{i}(h_i) \rightarrow \phi$ in the compact-open topology on $\Aut(K)$.
Since $\Inn(K)$ is a closed subgroup of $\Aut(K)$ it follows that $\phi$ is an
inner automorphism. This completes the proof of Lemmas \ref{centralizer} and
\ref{centralizer2}. \endproof

Let us return to the proof of \fullref{main}. Suppose we are in case (I).
Then $\gamma$ is contained in the compact subgroup $K K_1 \ldots K_r \le G$,
where $K_i = x_i K x_i^{-1}$.  Let $A$ be the closure in $G$ of the cyclic
group generated by $\gamma$. By our assumption, $A$ is compact.
Then $A\leq Z_{G}(\gamma )$, and
by \fullref{Juan}, if $w$ is a nontrivial element in $\Gamma _{2r}$
then $\{\alpha \in A \; : \; \rho _{\alpha }(w)\neq 1\}$
 is an open dense subset
of $A$. By Baire's theorem there is an $\alpha _{0}\in A$ such that $\rho
_{\alpha _{0}}$ is faithful. Its image contains $F$ and is isomorphic to the
surface group $\Gamma _{2r}.$ This completes the proof in case (I).

Finally suppose that we are in case (II).  Then $F^{\prime }\cap K=\{1\}$
and the sequence $(\pi \circ \rho _{\gamma ^{n}})_{n\geq 1}$ is eventually
faithful, where $\pi \co H\rightarrow H/K$ is the projection map. Let $\pi
(\gamma )=\beta (1)$ where $\{\beta (t)\}_{t}$ is a $1$-parameter subgroup
of $H/K$. The centralizer $Z_{H}(K)$ is closed in $H$, hence locally
compact, and by \fullref{centralizer}, $H/K\cong Z_{H}(K)/K\cap Z_{H}(K).$
\fullref{lift} ensures that $\{\beta (t)\}_{t}$ can be lifted to a $1$%
-parameter subgroup $\{c(t)\}_{t}$ in $Z_{H}(K)$, ie, $\beta (t)=c(t)K$. We
now claim that $\{c(t)\}_{t}\leq Z_{G}(\gamma )$.

Indeed, by \fullref{centralizer}, we can write $\gamma =ck$ where $c\in
Z_{H}(K)$ and $k\in K.$ However $\pi (\gamma )=\beta (1)$, hence $%
c(1)K=\gamma K,$ hence $c(1)k^{\prime }=c$ for some $k^{\prime }\in K.$ But
for each $t,$ $c(t)$ commutes with $c(1)$ and with $k^{\prime },$ hence it
commutes with $c.$ As $c(t)\in Z_{H}(K)$ it must also commute with $k$,
hence with $\gamma $. This proves the claim.

As a consequence, we obtain a one parameter family of representations $\rho
_{t}\co \Gamma _{2r}\rightarrow G$ by setting $\rho _{t}:=\rho _{c(t)}.$ Again $%
\{t\in \mathbb{R} \; : \; $ $\pi \circ \rho _{t}(w)\neq 1\}$ is open because $%
t\mapsto \rho _{t}(w)$ is continuous from $\mathbb{R}$ to $G.$ By \fullref{Juan}, the sequence $\pi \circ \rho _{n}=\pi \circ \rho _{\gamma ^{n}}$
is eventually faithful. This implies that the analytic map $t\mapsto \pi
\circ \rho _{t}(w)$ from $\mathbb{R}$ to the Lie group $H/K$ is not
constant. Therefore the set $\{t\in \mathbb{R} \; : \; \pi \circ \rho _{t}(w)\neq
1\}$ is dense. Again, by Baire's theorem there must be a $t_{0}\in \mathbb{R}
$ such that $\pi \circ \rho _{t_{0}}$ is faithful. Then $\rho _{t_{0}}$ is
also faithful and its image is a subgroup of $G$ isomorphic to $\Gamma _{2r}$
containing $F$. This completes the proof of \fullref{main}.




%

\section{The analytic structure of Hom$(\Gamma ,G)$}

In this section we will recall some facts about the structure of $%
\Hom(\Gamma,G)$ as an analytic variety, where $\Gamma$ is a finitely
generated group and $G$ is a Lie group.

Consider first the case that $\Gamma$ is isomorphic to a free group ${F}_k$
with free basis $e_1,\dots,e_k$. A homomorphism $\sigma\in\Hom({F}_k,G)$ is
determined by $\sigma(e_1),\dots,\sigma(e_k)$ and hence we have an
identification of $\Hom({F}_k,G)$ with the analytic manifold $%
G^k=G\times\dots\times G$. Given an element $\gamma=e_{i_1}\dots e_{i_l}\in{F%
}_k$ we consider the analytic map
\begin{equation*}
P_\gamma\co G^k\to G,\ \ P_\gamma(A_1,\dots,A_k)=A_{i_1}\cdots A_{i_l}.
\end{equation*}
The set $P_\gamma^{-1}(1_G)=\{\rho\in\Hom({F}_k,G) \; : \;  \gamma\in\Ker(\rho)\}$
is a closed analytic subvariety of $\Hom({F}_k,G)$. Recall the following
basic result (Epstein \cite{Eps}):

\begin{thm}
\label{free-generic} Let $G$ be a connected nonsolvable Lie group. Then
the set of faithful homomorphisms $\sigma \co {F}_{k}\to G$ is dense and has
full Haar measure in $\Hom({F}_k,G)\cong G^{k}$.
\end{thm}

We include a proof for the convenience of the reader.

\proof
By definition $\cup _{\gamma \in F_{k}\setminus 1_{F_{k}}}P_{\gamma }^{-1}(%
1_{G})$ is the complement of the set of faithful representations. The
claim follows from the Baire category theorem if $P_{\gamma }^{-1}(1_{G})$
is nowhere dense and of $0$ measure for all $\gamma \neq 1_{\Gamma }$.
Since $P_{\gamma }^{-1}(1_{G})$ is an analytic subvariety of $G^{k}$, it
is either nowhere dense and of $0$ measure or contains $G^{k}$ since $G$ is
connected. To show that the later case cannot occur, it suffices to find one
faithful representation, ie, it suffices to find a nonabelian free subgroup
of $G$. The existence of a free subgroup in $G$ follows from the
Tits alternative \cite{Tit}, or more simply from the fact that $G$ contains
a subgroup locally isomorphic to either $\text{PSL}_2({\mathbb{R}} )$ or 
$\text{PSO} (3)$, and each of these groups contains a free subgroup.
\endproof

Now let $\Gamma $ be a general finitely generated group. To a given
presentation $\Gamma =\langle \gamma _{1},\dots ,\gamma _{k}\ |\
\{R_{i}\}_{i\in I}\rangle $ of $\Gamma $, we associate the surjection $\pi \co {%
F}_{k}\to \Gamma $ defined by $\pi (e_{j})=\gamma _{j}$. The homomorphism $%
\pi $ induces an injective map
\begin{equation*}
\pi ^{*}\co \Hom(\Gamma,G)\to \Hom({F}_k,G)
\end{equation*}
and its image coincides with $\cap _{i\in I}P_{R_{i}}^{-1}(1_{G})$. Hence,
we can identify $\Hom(\Gamma,G)$ with an analytic subvariety of $G^{k}$. (In
fact, the induced structure of $\Hom(\Gamma,G)$ as an analytic variety does
not depend on the presentation of $\Gamma $. We will not use this.) An
important observation is that for all $\gamma \in \Gamma $ the map $%
P_{\gamma }^{\Gamma }\co \Hom(\Gamma,G)\to G$ given by $P_{\gamma }^{\Gamma
}(\rho )=\rho (\gamma )$ is analytic. Moreover, if $[\gamma ]\in {F}_{k}$ is
an element representing $\gamma $ then we have $P_{\gamma }^{\Gamma
}=P_{[\gamma ]}|_{\Hom(\Gamma,G)}$. This is why in the sequel we will
simplify notation and write $P_{\gamma }^{\Gamma }=P_{\gamma }$.

An important fact for our considerations is that analytic subvarieties admit
locally finite stratifications with smooth strata. The following crucial
result is due to Whitney, Thom and Lojasiewicz. We refer to Kaloshin \cite
{Kal} for its proof.

\begin{prop}
\label{prop:strata} Let $V$ be an analytic subvariety of an analytic
manifold $M$. Then there is a locally finite decomposition $V=\cup V_{i}$,
where $V_{i}$ are connected analytic submanifolds of $M$.
\end{prop}

The statement of \fullref{prop:strata} is much weaker than, and
follows directly from, \cite[Theorem 1]{Kal}. We have chosen this simplified
statement to avoid recalling the more subtle properties of stratifications.
Despite this, we will refer to the submanifolds $V_{i}$ as the \emph{strata}
of $V$.

%
%

\section{Dense subgroups of connected Lie groups}\label{labelled section}


This section will establish some properties of dense subgroups of connected
Lie groups. We begin by recalling some results from \cite{BG1}. We then
determine the number $d(G)$ of \fullref{frf}, and finally study the
structure of the set $\mathcal{D}(F_k ,G)$ of dense representations of the
free group $F_k$ in $G$.



A Lie group $H$ is \emph{topologically perfect} if its commutator group is
dense. Recall the following theorem from \cite{BG1} (see also \cite{GZ}):

\begin{thm}
\label{dense in TP} Let $H$ be a connected topologically perfect Lie group.
Assume that the Lie algebra $\text{Lie}(H)$ is generated (as a Lie algebra)
by $d=d(H)$ elements. Then there is an identity neighborhood $U\subset H$,
and a proper analytic subvariety $R\subset U^{d}$, such that $%
\langle h_{1},\ldots ,h_{d}\rangle $ is dense in $H$ for any $(h_{1},\ldots
,h_{d})\in U^{d}\setminus R$.
\end{thm}

When $G$ is topologically perfect we can define the constant $d(G)$ of
\fullref{frf} to be the minimal number of generators for $\text{Lie}(G)$%
. As a consequence of \fullref{dense in TP} we obtain:

\begin{cor} \label{dense in TP corollary}
Let $G$ be a connected topologically perfect Lie group. Then $\mathcal{D}%
(\Gamma ,G)$ is open in $\Hom (\Gamma ,G)$, and dense in a neighborhood of 
the trivial representation.
\end{cor}

\begin{proof}
If $\rho_0\in\Hom (\gC,G)$ is a representation of $\gC$ in $G$
with dense image, then for some $\gc_1,\ldots,\gc_d\in\gC$, where
$d=d(G)$, we have $(\rho_0(\gc_1),\ldots,\rho_0(\gc_d))\in
U^d\setminus R$. But then $(\rho (\gc_1),\ldots,\rho (\gc_d))\in
U^d\setminus R$ for any $\rho$ sufficiently close to $\rho_0$ in
$\Hom (\gC,G)$. By \fullref{dense in TP} any such $\rho$ has a
dense image.
\end{proof}

We now define $d(G)$ for a general connected Lie group $G$.

For a connected abelian Lie group $A,$ define $\text{rank}(A)$ as the
dimension of the tensor with ${\mathbb{R}}$, ie, if ${\mathbb{T}}$ is the
one dimensional torus and $A={\mathbb{T}}^{j}\times {\mathbb{R}}^{k}$ then $%
\text{rank}(A)=k$. It is easy to see that a generic set (in the sense of
the Baire category theorem or measure theory) of $k+1$ elements generates a
dense subgroup in $A$. For example, if A is compact then a generic element generates
a dense cyclic subgroup. In general, $k+1$ elements $a_{1},\ldots
,a_{k+1}\in A$ generate a dense subgroup if and only if the projections of
the first $k$ elements $a_{1},\ldots ,a_{k}$ to the second factor ${%
\mathbb{R}}^{k}$ form a basis, and, after identifying the compact quotient $%
A/\langle a_{1},\ldots ,a_{k}\rangle $ with ${\mathbb{T}}^{j+k}$, the $j+k$
coordinates of the projection of the last element $a_{k+1}$ to this torus
are independent. For a connected abelian Lie group $A$, we thus define $d(A)=%
\text{rank}(A)+1$.

Let now $G$ be a general connected Lie group. Set $G_0=G$ and define
inductively $G_i$ to be the closure of the derived group $[G_{i-1},G_{i-1}]$%
. (In more standard notation, the subgroup $G_i$ is denoted by the mildly
 cumbersome $\overline{G^{(i)}}$, which will not be	used here.)  
The decreasing sequence $G_i$ must stabilize after finitely many steps $m$
to a group $H=G_m$, and $H$ has the property that its commutator is dense,
ie, it is topologically perfect. The general case is reduced to the abelian
and topologically perfect cases using the following:

\begin{lem}
\label{123} A subgroup $D$ of $G$ is dense if and only if

\begin{enumerate}
\item  its image in $G/G_{2}$ is dense in $G/G_{2}$, and

\item  its intersection with $H$ is dense in $H$.
\end{enumerate}
\end{lem}

\begin{proof}
If $D$ is dense then (1) follows immediately. Moreover the
commutator group $[D,D]$ is clearly dense in $G_1$, and by a
simple induction the $m^{\text{th}}$ commutator of $D$ is dense in
$H=G_m$.

The other direction will follow if we can show that (1) implies 
the image of $G$ in $G/H$ is dense in $G/H$. To do this we will use
the fact that $G/H$ is solvable.
In a connected solvable Lie group $B$, a subgroup is dense if and only if
its image in $B/B_2$ (modulo the second closed commutator) is
dense in $B/B_2$. To see this, note that the commutator of a
connected solvable Lie group is nilpotent, and that a subgroup of
a nilpotent group is dense if and only if it is dense modulo the first
commutator.
\end{proof}

We define the number $d(G)$ as follows\footnote{%
In case $G/H$ is nilpotent, one can take $d(G)=\max \{d(G_{0}/G_{1}),d(H)\}.$%
}:
\begin{equation*}\label{defdg}
d(G)=\max \{d(G_{0}/G_{1}),d(G_{1}/G_{2}),d(H)\},
\end{equation*}
where $d(H)$ is the minimal number of generators for the Lie algebra of $H$,
and $d(G_{i}/G_{i+1})=\text{rank}(G_{i}/G_{i+1})+1$.

Consider the following subsets of $\Hom (\Gamma,G)$:
\begin{equation*}
\mathcal{D}_H(\Gamma,G)=\{\rho\in\Hom (\Gamma,G) \; : \; \overline{\rho
(\Gamma)\cap H}=H\},~\text{and}
\end{equation*}
\begin{equation*}
\mathcal{D}_{G/G_2}(\Gamma,G)=\{\rho\in\Hom (\Gamma,G) \; : \; \overline{\rho
(\Gamma)G_2}=G\}.
\end{equation*}
By \fullref{123} we have:
\begin{equation*}
\mathcal{D}(\Gamma,G)=\mathcal{D}_H(\Gamma,G)\cap\mathcal{D}%
_{G/G_2}(\Gamma,G).
\end{equation*}
Moreover, for free groups we have:

\begin{lem}
\label{6.6} Suppose that $k\geq d(G)$, then:

\begin{itemize}
\item  The set $\mathcal{D}_{H}({F}_{k},G)$ is open in $\Hom
({F}_{k},G)$.

\item  The set $\mathcal{D}_{G/G_{2}}({F}_{k},G)$ is the complement of a
countable union of proper closed analytic subvarieties of $\Hom ({F}_{k},G)$%
. In particular, it is of second category.
\end{itemize}
\end{lem}

\proof
The first claim follows from \fullref{dense in TP}: if $\rho _{0}\in
\mathcal{D}_{H}({F}_{k},G),$ then $\rho (F_{k})_{m}$ is dense in $G_{m}=H.$
Hence there are $d(H)$ words involving commutators of length $m$ in $k$
letters, such that when applying them to the image (under $\rho _{0}$) of
the generators of ${F}_{k}$ (think of them as the coordinates of a point
in $G^{d(H)}$) one gets a point in $U^{d(H)}\setminus R\subset H^{d(H)}$.
Clearly if $\rho $ is sufficiently close to $\rho _{0}$ then the same words
applied to the $\rho $ image of the generators still yield a point in $%
U^{d(H)}\setminus R$. By \fullref{dense in TP}, $\rho ({F}_{k},G)\cap H$
is dense in $H$.

To see the second claim note that a subgroup of $G/G_2$ is dense in $G/G_2$
if its image in $G/G_1$ is dense in $G/G_1$ and its intersection with $G_1$
projects to a dense subgroup of $G_1/G_2$. Both conditions are generic in
the sense that their complements are a countable union of proper analytic
closed subvarieties: there are $d(G_0/G_1)$ words with $k$ letters which
generically generate a dense subgroup in the quotient $G_0/G_1$, and $%
d(G_1/G_2)$ words involving commutators of the $k$ letters which
generically generate a dense subgroup in the quotient $G_1/G_2$.
The former assertion is clear, while the latter is a little harder
to see and we leave it to the reader as an exercise. In fact, if
$x_1, ..., x_k$ are generic elements of G then the commutators
$[x_1,x_i]$ for $2\leq i \leq k$ form a basis of $G_1/G_2$, and
together with $\prod_{i=2}^k[x_1^i,x_i]$ they generate a dense
subgroup of $G_1/G_2$.
\endproof

\begin{cor}\label{6.6 corollary}
Assume $G$ is a connected nonsolvable Lie group. If $k \ge d(G)$ then there exists a sequence of faithful representations 
$(\sigma_i) \subset \mathcal{D} (F_k,G)$ converging to the trivial representation.
\end{cor}

\proof
By \fullref{dense in TP corollary} there exists a sequence $(\sigma'_i) 
\subset \mathcal{D} (F_k,H) $ converging to the trivial representation.
By \fullref{6.6} and \fullref{free-generic}, it is possible to obtain the 
desired sequence $( \sigma_i)$ via an arbitrarily small perturbation of
$(\sigma'_i)$.
\endproof




\section[Proof of \ref{frf}]{Proof of \fullref{frf}}

We are now in a position to complete the proof of \fullref{frf}. We
begin by fixing once and for all a relatively compact open neighborhood $%
B\subset \Hom(\Gamma,G)$ of the trivial homomorphism and let $V_{1},\dots
,V_{s}$ be finitely many strata of $\Hom(\Gamma,G)$ covering $B$.

The group $\Gamma$ is, as assumed, $d$-fully residually free for $d=d(G)$.
Hence it is generated by $k\geq d$ elements and there is a sequence of
surjective homomorphisms $\phi_i\co \Gamma\to{F}_d$ such that for every $%
\gamma\in\Gamma\setminus1_\Gamma$ there is $i_\gamma$ with $\gamma\notin%
\Ker(\phi_i)$ for all $i\ge i_\gamma$. The assumption that $\Gamma$ is
nonabelian implies that $d\ge 2$. The homomorphisms $\phi_i$ induce analytic
maps $\phi_i^*\co \Hom({F}_d,G)\to\Hom(\Gamma,G)$.  By \fullref{6.6 corollary}
we can choose a sequence $(\sigma_i)\subset\Hom({F}_d,G)$
 of faithful representations with dense image
sufficiently close enough to the trivial homomorphism so that $%
\phi_i^*(\sigma_i)\in B$ for all $i$. Up to passing to a subsequence and
relabelling, we may assume that $\phi_i^*(\sigma_i)\in V_1\subset\Hom(\gC,G)$
for all $i$.

Given $\gamma \in \Gamma \setminus 1_{\Gamma }$ we deduce from the
connectivity of $V_{1}$, using analytic continuation and the implicit
functions theorem, that either $V_{1}\subset \{P_{\gamma }^{-1}(1_{G})\}$
or $V_{1}\cap \{P_{\gamma }^{-1}(1_{G})\}$ is nowhere dense in $V_{1}$.
The former case cannot occur since by construction we have $P_{\gamma }(\phi
_{i}^{*}(\sigma _{i}))=\sigma _{i}(\phi _{i}(\gamma ))\neq 1_{G}$ for all $%
i$ sufficiently large. In particular, Baire's category theorem implies that
the set $\mathcal{F}(\Gamma ,G)\cap V_{1}$ of all faithful $\rho \in V_{1}$
is of second category in $V_{1}$.

By construction the image of $\phi_1^*(\sigma_1)$ coincides with the image
of $\sigma_1$ and hence is dense. This implies that the open subset $%
\mathcal{D}_H(\Gamma ,G)\cap V_1$ of $V_1$ is nonempty.
Additionally it implies that the set $\mathcal{D}_{G/G_{2}}(\Gamma ,G)\cap
V_{1}$ is nonempty.  Since it is the complement of a countable union of
proper closed analytic subvarieties we conclude, again by analyticity and
the implicit functions theorem, that all these varieties are proper, and
that $\mathcal{D}_{G/G_{2}}(\Gamma ,G)\cap V_{1}$ is also of second
category. Hence the intersection
\begin{equation*}
\mathcal{F}(\Gamma ,G)\cap \mathcal{D}(\Gamma ,G)\cap V_{1}=\mathcal{F}%
(\Gamma ,G)\cap \mathcal{D}_{G/G_{2}}(\Gamma ,G)\cap \mathcal{D}_{H}(\Gamma
,G)\cap V_{1}
\end{equation*}
is not empty because it is the intersection of the second category subset
\begin{equation*}
\mathcal{F}(\Gamma ,G)\cap \mathcal{D}_{G/G_{2}}(\Gamma ,G)\cap V_{1}
\end{equation*}
with the nonempty open subset $\mathcal{D}_{H}(\Gamma ,G)\cap V_{1}$ of $V_{1}$%
. \qed




\section{Some remarks on connected compact groups}

When $G$ is a connected compact Lie group then $d(G)=1$ if $G$ is abelian
and $d(G)=2$ if it is not. For compact semisimple Lie groups one can easily
deduce the following lemma from \fullref{dense in TP}. The general case
follows by an simple argument similar to the one given in the proof of \fullref{6.6}.

\begin{lemma}
\label{css} Let $G$ be a connected compact Lie group. Then the set $\mathcal{%
D}({F}_{2},G)$ is of full Haar measure and Baire dense in $%
\Hom({F}_2,G)\cong G\times G.$
\end{lemma}

When $G$ is nonabelian, \fullref{free-generic} says that a generic pair $%
(a,b)\in G\times G$ also generates a free group. We shall now generalize
this result to an arbitrary compact connected group.

\begin{proposition}
\label{fscg}Let $G$ be a second countable connected compact nonabelian
group. Then there exists a subset $\mathcal{O}$ in $G\times G,$ which is
both of second Baire category and of full Haar measure, such that any pair $%
(a,b)$ in $\mathcal{O}$ generates a dense free subgroup in $G$.
\end{proposition}

\begin{proof}
By the Peter-Weyl theorem (c.f. \cite{MZ}) there is a decreasing
sequence of normal compact subgroups $K_{n}\lhd G$ such that each
quotient $G/K_{n}$ is a connected compact Lie group and
$\bigcap_{n\geq 1}K_{n}=\{1\}$. Let $\mathcal{D}_{n}$ be the set
of all pairs $(a,b)$ in $G\times G$ that generate a dense subgroup
in the quotient $G/K_{n}.$ Clearly $\mathcal{D}({F}_2,G)=\bigcap
\mathcal{D}_{n}$, and by \fullref{css}, $\mathcal{D}_{n}$ is of
second Baire category and of full Haar measure.

The analogous assertion for $\mathcal{F}({F}_2,G)$ follows easily from \fullref{free-generic} because one of the quotients $G/K_n$ is
nonabelian.
\end{proof}

\begin{cor}
Any connected second countable nonabelian compact group contains a dense
surface group of genus 2.
\end{cor}

\bibliographystyle{gtart}
\bibliography{link}

\end{document}